\documentclass[twoside]{article}
\usepackage{amsfonts,amsmath,amsthm}
\usepackage[latin1]{inputenc}

\def\qed{\hbox to 0pt{}\hfill$\rlap{$\sqcap$}\sqcup$}

\newtheorem{thm}{Theorem}
\usepackage{epsfig}

\makeatletter

\begin{document}

\title{An extension of the characterization of the domain of attraction of an asymptotically  stable
fixed point in the case of a nonlinear discrete dynamical system}
\author{E. Kaslik$^{1,3}$, A.M. Balint$^{2}$, A. Grigis $^{3}$, St. Balint$^{1}$ }

\date{}
\maketitle

\noindent$^1$Department of Mathematics, West University of
Timi\c{s}oara, Romania \hfill\break e-mail: balint@balint.uvt.ro

\noindent$^2$Department of Physics, West University of
Timi\c{s}oara, Romania

\noindent$^3$L.A.G.A, Institut Galilee, Universite Paris 13,
France

\begin{abstract}
The purpose of this paper is to show that if the spectral radius of the matrix of the
linearized system of a nonlinear discrete dynamical system is less then one then the characterization theorem from \cite{KBBB1,KBGB2} holds.
\end{abstract}



\section{Introduction}

We consider the system of difference equations
\begin{equation}
\label{ec} x_{k+1}=g(x_{k})\qquad k=0,1,2...
\end{equation}
where $g:\Omega\rightarrow\Omega$ is an analytic function defined
on a domain $\Omega$ included in $\mathbb{R}^{n}$.

A point $x^{0}\in\Omega$ is a fixed point for the system
(\ref{ec}) if $x^{0}$ satisfies
\begin{equation}\label{fix}
    x^{0}=g(x^{0})
\end{equation}

The fixed point $x^{0}$ of (\ref{ec}) is "stable" provided that
given any ball
$B(x^{0},\varepsilon)=\{x\in\Omega/\|x-x^{0}\|<\varepsilon\}$,
there is a ball
$B(x^{0},\delta)=\{x\in\Omega/\|x-x^{0}\|<\delta\}$ such that if
$x\in B(x^{0},\delta)$ then $g^{k}(x)\in B(x^{0},\varepsilon)$,
for $k=0,1,2,...$ \cite{Kelley-Peterson}.

If in addition there is a ball $B(x^{0},r)$ such that
$g^{k}(x)\rightarrow x^{0}$ as $k\rightarrow\infty$ for all $x\in
B(x^{0},r)$ then the fixed point $x^{0}$ is "asymptotically
stable".\cite{Kelley-Peterson}.

The domain of attraction $DA(x^{0})$ of the asymptotically stable
fixed point $x^{0}$ is the set of initial states $x\in \Omega$
from which the system converges to the fixed point itself i.e.
\begin{equation}\label{da}
    DA(x^{0})=\{x\in\Omega | g^{k}(x)\stackrel{k\rightarrow\infty}{\longrightarrow}x^{0}\}
\end{equation}

It is known that $x^{0}$ is a fixed point for system (\ref{ec}) if
and only if $0\in\mathbb{R}^{n}$ is a fixed point for the system
\begin{equation}
\label{ec0} y_{k+1}=f(y_{k}) \qquad k=0,1,2...
\end{equation}
where $f:\Omega-x^{0}\rightarrow\Omega-x^{0}$ is the analytic
function defined by
\begin{equation}\label{f}
    f(y)=g(y+x^{0})-x^{0} \qquad \textrm{for } y\in\Omega-x^{0}
\end{equation}

The fixed point $x^{0}$ of (\ref{ec}) is asymptotically stable if
and only if the fixed point $0\in\mathbb{R}^{n}$ of the system
(\ref{ec0}) is asymptotically stable.

The domain of attraction of $x^{0}$, $DA(x^{0})$ is related to the
domain of attraction of $0$, $DA(0)$ by the equation
\begin{equation}\label{relda}
    DA(x^{0})=DA(0)+x^{0}
\end{equation}

For the above reason in the followings instead of the system
(\ref{ec}) we will consider the system (\ref{ec0}).

Theoretical research shows that the $DA(0)$ and its boundary are
complicated sets. In most cases, they do not admit an explicit
elementary representation. For this reason, different procedures
are used for the approximation of the $DA(0)$ with domains having
a simpler shape. For example, in the case of the theorem 4.20 pg.
170 \cite{Kelley-Peterson} the domain which approximates the
$DA(0)$ is defined by a Lyapunov function $V$ built with the
matrix $\partial_{0}f$ of the linearized system in $0$. In
\cite{KBGB2}, we presented a technique for the construction of a
Lyapunov function $V$ in the case when the matrix $\partial_{0}f$
is a contraction, i.e. $\|\partial_{0}f\|<1$. In this paper, we
extend this result for the more general case when
$r(\partial_{0}f)<1$ (where $r$ denotes the spectral
radius). The Lyapunov function $V$ is built using the whole
nonlinear system, not only the matrix $\partial_{0}f$. $V$ is
defined on the whole $DA(0)$, and more, the $DA(0)$ is the natural
domain of analyticity of $V$.

\section{The result}

Let be $f:\Omega\rightarrow\Omega$ an analytic function defined on
a domain $\Omega\subset\mathbb{R}^{n}$ containing the origin
$0\in\mathbb{R}^{n}$.

\begin{thm}
If the function $f$ satisfies the following conditions:
\begin{equation}
  f(0) = 0
\end{equation}
\begin{equation}
  r(\partial_{0}f) < 1
\end{equation}
then $0$ is an asymptotically stable fixed point. $DA(0)$ is an
open subset of $\Omega$ and coincides with the natural domain of
analyticity of the unique solution $V$ of the iterative first
order functional equation
\begin{equation}
\label{ecV}
\begin{array}{ll}
\left\{\begin{array}{l}
V(f(x))-V(x)=-\|x\|^{2}\\
V(0)=0
\end{array}\right.
\end{array}
\end{equation}
The function $V$ is positive on $DA(0)$ and
$V(x)\stackrel{x\rightarrow x^{0}}{\longrightarrow}+\infty$, for
any $x^{0}\in FrDA(0)$ ($FrDA(0)$ denotes the boundary of
$DA(0)$).
\end{thm}

\begin{proof}Let be $A=\partial_{0}f$ and
$r=r(\partial_{0}f)$.

The fact $0$ is an asymptotically stable fixed point is proved in
\cite{Elaydi}. The fact that $DA(0)$ is an open subset of $\Omega$ follows
from the continuity of the function $f^{k}$, for any $k\in\mathbb{N}$.

We will prove that the series $\sum\limits_{k=0}^{\infty}\|f^{k}(x)\|^{2}$ is convergent for any
$x\in DA(0)$.

We can decompose the function $f$ as follows:
\begin{equation}
    f(x)=Ax+h(x)
\end{equation}
where $h:\Omega\rightarrow\Omega$ is an analytical function which
satisfies $h(0)=0$.

As $r<1$, there exists $c>0$ such that
\begin{equation}\label{cond A}
    \|A^{k}x\|<cr^{k}\|x\|\qquad \forall x\in\Omega, k\in\mathbb{N}
\end{equation}

Let be $\bar{c}=\max\{1,c\}$ and $\varepsilon=\frac{1-r}{2\bar{c}}>0$. As $h$ is continuous and $h(0)=0$, it follows that
there exists $\delta>0$ such that
\begin{equation}\label{cond h}
    \|h(x)\|<\varepsilon\|x\| \qquad \forall x\in B(0,\delta).
\end{equation}

Let be $x\in DA(0)$ and $x_{k}=f^{k}(x)$. This provides the existence of $k_{x}\in\mathbb{N}$ such that
$x_{k}\in B(0,\delta)$, for any $k\geq k_{x}$. Let be $y_{k}=x_{k+k_{x}}$. This sequence satisfies
$y_{k}\in B(0,\delta)$ for any $k\in\mathbb{N}$ and it also satisfies (\ref{ec0}).

The formula of variation of constants gives:
\begin{equation}\label{voc}
    y_{k}=A^{k}y_{0}+\sum_{i=l}^{k-1}A^{k-i-1}h(y_{i})\qquad \forall k\in\mathbb{N^{\star}}
\end{equation}

Relations (\ref{cond A}) and (\ref{voc}) provide
\begin{equation}\label{ineg1}
    \|y_{k}\|\leq cr^{k}\|y_{0}\|+\sum_{i=0}^{k-1}cr^{k-i-1}\|h(y_{i})\| \qquad \forall k\in\mathbb{N^{\star}}
\end{equation}
and using (\ref{cond h}) and that $c\leq \bar{c}$, the following inequality follows:
\begin{equation}\label{ineg2}
    \|y_{k}\|\leq \bar{c}r^{k}\|y_{0}\|+\sum_{i=0}^{k-1}\bar{c}r^{k-i-1}\varepsilon\|y_{i}\|\qquad \forall k\in\mathbb{N^{\star}}
\end{equation}

Relation (\ref{ineg2}) can be written as
\begin{equation}\label{ineg2'}
    r^{-k}\|y_{k}\|\leq \bar{c}\|y_{0}\|+\sum_{i=0}^{k-1}\bar{c}r^{-1}\varepsilon(r^{-i}\|y_{i}\|) \qquad \forall k\in\mathbb{N^{\star}}
\end{equation}

Gronwall's inequality for the discrete case provides
\begin{equation}\label{ineg3}
    r^{-k}\|y_{k}\|\leq \bar{c}\|y_{0}\|\prod_{i=0}^{k-1}(1+\bar{c}r^{-1}\varepsilon)\qquad \forall k\in\mathbb{N^{\star}}
\end{equation}
thus
\begin{equation}\label{ineg3'}
    \|y_{k}\|\leq r^{k}\bar{c}\|y_{0}\|(1+\bar{c}r^{-1}\varepsilon)^{k}=
    \bar{c}\|y_{0}\|(r+\bar{c}\varepsilon)^{k}=\bar{c}\|y_{0}\|(\frac{r+1}{2})^{k} \qquad \forall k\in\mathbb{N^{\star}}
\end{equation}

Denoting $\frac{r+1}{2}=\alpha<1$, the relation (\ref{ineg3'}) gives
\begin{equation}\label{ineg}
    \|x_{k}\|\leq\bar{c}\|x_{k_{x}}\|\alpha^{k}\qquad \forall k\geq k_{x}
\end{equation}

Thus, for any $x\in DA(0)$ there exists $k_{x}\geq 0$ such that
\begin{equation}\label{ineg'}
    \|f^{k}(x)\|\leq\bar{c}\|f^{k_{x}}(x)\|\alpha^{k}\qquad \forall k\geq k_{x}
\end{equation}
which assures that the series $\sum\limits_{k=0}^{\infty}\|f^{k}(x)\|^{2}$ is convergent for any
$x\in DA(0)$.

Let be $V=V(x)$ the function defined by
\begin{equation}\label{t3}
    V(x)=\sum\limits_{k=0}^{\infty}\|f^{k}(x)\|^{2} \qquad \forall
    x\in DA(0)
\end{equation}

The above function defined on $DA(0)$ is analytical, positive and
satisfies (\ref{ecV}). In order to show that the function $V$
defined by (\ref{t3}) is the unique function which satisfies
(\ref{ecV}) we consider $V'=V'(x)$ satisfying (\ref{ecV}) and we
denote by $V''$ the difference $V''=V-V'$. It is easy to see that
$V''(f(x))-V''(x)=0$, for any $x\in DA(0)$. Therefore, we have
$V''(x)=V''(f^{k}(x))$ for any $x\in DA(0)$ and any $k\in\mathbb{N}$.
It follows that
$V''(x)=\lim\limits_{k\rightarrow\infty}V''(f^{k}(x))=0$ for any
$x\in DA(0)$. In other words, $V(x)=V'(x)$, for any $x\in DA(0)$,
so $V$ defined by (\ref{t3}) is the unique function which
satisfies (\ref{ecV}).

In order to show that $V(x)\stackrel{x\rightarrow
x^{0}}{\longrightarrow}\infty$ for any $x^{0}\in  FrDA(0)$ we
consider $x^{0}\in  FrDA(0)$ and $r>0$ such that
$\|f^{k}(x^{0})\|>r$, for any $k=0,1,2...$. For an arbitrary
positive number $N>0$ we consider the first natural number $k_{1}$
which satisfies $k_{1}\geq\frac{2N}{r^{2}}+1$. Let be $r_{1}>0$
such that $\|f^{k}(x)\|\geq\frac{r}{\sqrt{2}}$ for any
$k=1,2,..,k_{1}$ and $x\in B(x^{0},r_{1})$. For any $x\in
B(x^{0},r_{1})\cap DA(0)$ we have
$\sum\limits_{k=0}^{k_{1}}\|f^{k}(x)\|^{2}>N$. Therefore,
$V(x)\stackrel{x\rightarrow x^{0}}{\longrightarrow}\infty$.
\end{proof}

\section{Numerical examples}
In this section, some examples of discrete systems are given, for which the domain of
attraction of the zero steady state can be estimated using the method described in \cite{KBBB1}.

\subsubsection*{Example 1}

The following system of difference equations is
considered:
\begin{equation}
\label{ex3}
\begin{array}{ll}
\left\{\begin{array}{l}
x_{n+1}=x_{n}y_{n}+y_{n}\\
y_{n+1}=y_{n}^{3}
\end{array}\right.
\end{array}
\end{equation}
\noindent It is clear that $(0,0)$ is an asymptotically stable
fixed point for the system (\ref{ex3}). It can be proved theoretically that the domain of attraction of $(0,0)$ is
$DA(0)=\mathbb{R}\times (-1,1)$.

\noindent The matrix of the linearized system
has the norm $\|\partial_{0}f\|=1$, thus the theorem from \cite{KBBB1,KBGB2} does not apply. On the other
hand, $r(\partial_{0}f)=0$, thus the characterization theorem for $DA(0)$ given above applies. Thus,
the domain of attraction can be estimated using the numerical method given in \cite{KBBB1}. The first step of the
method gives the whole domain of attraction:
\renewcommand{\thefigure}{\arabic{figure}}
\begin{figure}[htbp]
\centering
\includegraphics*[bb=3cm 0cm 13.5cm
10.5cm, width=5cm]{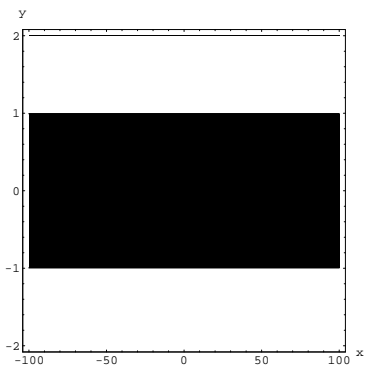} \caption{The estimate of $DA(0)$ after
1 step for system (\ref{ex3})}
\end{figure}

\subsubsection*{Example 2}
The following system is considered:
\begin{equation}
\label{ex5}
\begin{array}{ll}
\left\{\begin{array}{l}
x_{n+1}=4x_{n}^{2}+y_{n}\\
y_{n+1}=x_{n}y_{n}
\end{array}\right.
\end{array}
\end{equation}
\noindent This system of difference equations has the asymptotically
stable fixed point $(0,0)$. It is once again the case when $\|\partial_{0}f\|=1$ and $r(\partial_{0}f)=0$.
After one step, we can obtain the following estimate of the $DA(0)$:
\begin{figure}[htbp]
\centering
\includegraphics*[bb=3cm 0cm 13.5cm
10.5cm, width=5cm]{fig2D_RS'.eps} \caption{The estimate of $DA(0)$ after
1 step for system (\ref{ex5})}
\end{figure}

\bibliography{Evabib}
\end{document}